\input ./style/arxiv-general.cfg
\documentclass[MSNbibl,nameyear,seceqn,dvips]{arxstspdf}
\makeatletter
   \@ifpackageloaded{graphicx}{}{\usepackage{graphicx}}
\makeatother
\usepackage{flushend}
\usepackage{stfloats}

\volume{30}
\issue{2}
\pubyear{2015}
\firstpage{216}
\lastpage{227}
\doi{10.1214/14-STS507}
\docsubty{FLA}

\makeatletter
\newcommand{\ds}{\displaystyle}
\newcommand{\eqref}[1]{(\ref{#1})}
\newtheorem{theorem}{Theorem}[section]
\newtheorem{proposition}[theorem]{Proposition}
\newproclaim{remark}[theorem]{Remark}
\makeatother

\begin{document}
\begin{frontmatter}

\title{On Various Confidence Intervals Post-Model-Selection}%
\runtitle{On Various Confidence Intervals Post-Model-Selection}

\begin{aug}
\author[A]{\fnms{Hannes}~\snm{Leeb}\corref{}\ead[label=e1]{hannes.leeb@univie.ac.at}},
\author[B]{\fnms{Benedikt M.}~\snm{P\"otscher}\ead[label=e2]{benedikt.poetscher@univie.ac.at}}
\and
\author[C]{\fnms{Karl}~\snm{Ewald}\ead[label=e3]{karl.ewald@tuwien.ac.at}}
\runauthor{H. Leeb, B. M. P\"otscher and K. Ewald}

\affiliation{University of Vienna, University of Vienna and Vienna
University of Technology}

\address[A]{Hannes Leeb is Professor, Department of Statistics, University of Vienna,
Oskar-Morgenstern-Platz 1, 1090 Vienna, Austria \printead{e1}.}
\address[B]{Benedikt M. P\"otscher is Professor, Department of Statistics, University of Vienna, Oskar-Morgenstern-Platz 1, 1090 Vienna, Austria \printead{e2}.}
\address[C]{Karl Ewald is Ph.D. Candidate, Institute for Mathematical Methods in Economics,
Vienna University of Technology, Argentinierstrasse 8/E105-2, 1040 Vienna  \printead{e3}.}
\end{aug}

%
\begin{abstract}
We compare several confidence intervals after model selection
in the setting recently studied by
Berk et~al. [\textit{Ann. Statist.} \textbf{41} (2013) 802--837],
where the goal is to cover
not the true parameter but
a certain nonstandard quantity of interest
that depends on the selected model.
In particular, we compare the PoSI-intervals that are
proposed in that reference with the ``naive'' confidence interval,
which is constructed as if the selected model were correct and fixed a priori
(thus ignoring the presence of model selection).
Overall, we find that
the actual coverage probabilities of all these intervals deviate
only moderately from the desired nominal coverage probability.
This finding is in stark contrast to several papers in the existing literature,
where the goal is to cover the true parameter.
\end{abstract}

%
\end{frontmatter}

\section{Introduction and Overview}
\label{intro}

There is ample evidence in the literature that
model selection can have a detrimental impact on
subsequently constructed inference procedures like
confidence sets, if these are constructed in the ``naive'' way
where the presence of model selection is ignored.
Such results are reported, for example, by
\citet{Bro67a};
\citet{Bue63a};
\citet{Dij88a};
\citeauthor{Kab98a} (\citeyear{Kab98a},
\citeyear{Kab09a});
\citet{Kab06a};
\citet{Lee01d};
\citeauthor{Lee02d}
(\citeyear{Lee02d}, \citeyear{Lee03a},
\citeyear{Lee02a},
\citeyear{Lee02c},
\citeyear{Lee02b}, \citeyear{Lee06a});
\citet{Ols73a}; 
\citeauthor{Poe91a}
(\citeyear{Poe91a}, \citeyear{Poe06a});\break 
\citet{Poe09b};
\citeauthor{Poe09c}
(\citeyear{Poe09c}, \citeyear{Poe10a},
\citeyear{Poe11a});
\citet{Sen79a};
\citet{Sen87a}.

Recently, \citet{Ber13a} proposed a new class of confidence intervals,
so-called PoSI-intervals, which correct for the presence of model selection,
in the sense that these intervals guarantee a user-specified minimal coverage
probability, even if the model has been selected in a data-driven way.
However, the setting of \citet{Ber13a} differs from earlier studies,
in that they consider confidence intervals for
\textit{a different quantity of interest}:
In the aforementioned analyses, the quantity of interest
(the coverage target)
is always
a~fixed parameter or subparameter of the data-generating model.
In \citet{Ber13a}, on the other hand, a different and nonstandard
coverage target is considered
that depends on the selected model.
(Even if an overall correct model is assumed,
that nonstandard coverage target does \textit{not} coincide with
a parameter in the model, except for degenerate and trivial
situations.)
By design, the PoSI-intervals hence
do not provide a solution to the more traditional problem, where the
goal is to cover a parameter in the overall model after model selection.

\citet{Ber13a} motivate the need for PoSI-intervals by
the poor performance of the ``naive'' interval as observed in
the studies mentioned in the first paragraph of this section.
However, these studies do \textit{not} deal with the performance
of the ``naive'' procedures post-model-selection when the coverage target
is as in \citet{Ber13a}.
This raises the question of how the ``naive'' interval performs
when it is used to cover the coverage target considered in \citet{Ber13a}.\vadjust{\goodbreak}
The main contribution of this paper is to answer this.
In particular, we compare ``naive'' confidence intervals
and PoSI-intervals in the setting of \citet{Ber13a}.
[The results in the present paper are partly based on \citet{Ewa12a}, and
we refer to this thesis for additional results and discussion.]

We find that the minimal coverage probability of the ``naive'' interval
is slightly below the nominal one, while that
of the various PoSI intervals is slightly above, when the coverage target
is as in \citet{Ber13a} and when AIC, BIC, or the LASSO are used
for model selection. In the scenarios that we consider, the coverage
probabilities of all these intervals are mostly within 10\% of the nominal
coverage probability.
In the more traditional setting where the coverage target is a
parameter in the overall model, however,
all these intervals generally fail to deliver the desired
minimal coverage probability. (Note that the various
PoSI-intervals are not designed to deal with this coverage target.)
For example, consider the scenario depicted by the solid
curves in Figure~\ref{fig1} on page~\pageref{fig1}:
There, a~``naive'' confidence interval post-model-selection
with nominal coverage probability 0.95
has a minimal coverage probability of about 0.91 and
the corresponding PoSI-interval has a minimal coverage probability
of about 0.96, if the coverage target is as in \citet{Ber13a}.
But if the coverage target is a parameter in the overall model,
the minimal coverage probabilities of the ``naive'' interval and
of the PoSI-interval drop to about 0.56 and 0.62, respectively.

The paper is organized as follows:
In Section~\ref{setting} we introduce the data-generating process,
the model-selection procedures,
the coverage targets, and
various confidence procedures, including the PoSI-intervals.
We consider the same assumptions
and constructions as in~\citet{Ber13a},
as well as some additional confidence intervals.
The (minimal) coverage probabilities of ``naive'' intervals
and of PoSI-intervals are studied in Sections~\ref{theory}
and~\ref{simulation}.
In particular, Section~\ref{theory} contains an explicit finite-sample
analysis of these procedures in a simple scenario
with two nested candidate models.
Section~\ref{simulation} contains a simulation study where we compare these
intervals in three more complex scenarios; the first scenario is also studied
by \citet{Kab06a}, and the other two scenarios are taken from \citet{Ber13a}.
(The code used for the computations in Section~\ref{theory} and for
the simulations in Section~\ref{simulation} is available from the first
author on request.)
Finally, in the \hyperref[app]{Appendix} we present an example
with a coverage target that is similar to, but slightly different
from, that considered in \citet{Ber13a}. The interesting
feature of this example
is that the ``naive''\vadjust{\goodbreak} confidence interval here is valid,
in the sense that its coverage probability is
never below the nominal level.

\section{Coverage Targets and Confidence~Intervals}
\label{setting}

Throughout, we consider a set of $n$ homoskedastic Gaussian observations
with mean vector $\mu\in\mathbb{R}^n$ and common variance $\sigma
^2>0$, that is,
%
\begin{equation}
\label{model}
y = \mu+ u,
\end{equation}
where $u \sim N(0,\sigma^2 I_n)$.
We further assume that we have an estimator $\hat{\sigma}^2$ for
$\sigma^2$ that is independent of all the least-squares estimators that
will be introduced shortly.
See Remark~\ref{r2}(ii)
for some cautionary comments regarding our assumptions on $\hat{\sigma}^2$.
For the estimator $\hat{\sigma}^2$,
we either assume that it is distributed as a chi-squared random variable
with $r$ degrees of freedom multiplied by $\sigma^2 / r$, that is,
$\hat{\sigma}^2 \sim\sigma^2 \chi^2_r / r$, for some $r\geq1$,
or we assume that the variance is known a priori, in which case we set
$\hat{\sigma}^2 = \sigma^2$ and $r=\infty$.
Unless noted otherwise, all considerations that follow apply
to both the known-variance case and the unknown-variance case.
The joint distribution of $y$ and $\hat{\sigma}$ depends on
the parameters $\mu\in\mathbb{R}^n$ and $\sigma>0$,
and will be denoted by~$\mathbb{P}_{\mu,\sigma}$.

The available explanatory variables are represented by the columns of
a fixed
$n \times p$ matrix $X$, where we allow for $p>n$;
again, see Remark~\ref{r2}(ii).
We consider models where $y$ is regressed on a (nonempty) subset of
the regressors in $X$: For each model
$M \subseteq\{1,\ldots, p\}$ with $M \neq\varnothing$,
write $X_M$ for the matrix of those
columns of $X$ whose indices lie in $M$.
Writing $M$ as $M = \{ j_1,\ldots, j_{|M|}\} \subseteq\{1,\ldots, p\}$,
we thus have $X_M = (X_{j_1},\ldots, X_{j_{|M|}})$, where
$X_j$ denotes the $j$th column of $X$
and where $|M|$ denotes the size of $M$.
Write $\mathcal{M}$ for a user-specified (nonempty)
collection of candidate models.
Throughout, we assume that $\mathcal{M}$ consists only of
submodels of full column rank,
that is, we assume that the rank of $X_M$ equals $|M|$ and satisfies
$1 \leq|M| \leq n$ for each $M\in\mathcal{M}$.

Under a candidate model $M \in\mathcal{M}$, $y$ is modeled as
\[
y = X_M \beta_M + v_M,
\]
where $\beta_M$ corresponds to the orthogonal projection of~$\mu$
from \eqref{model}
onto the column-space of $X_M$, that is,
$\beta_M = (X_M' X_M)^{-1} X_M' \mu$.
The least-squares estimator corresponding to the model $M$ will be denoted
by $\hat{\beta}_M$, that is, $\hat{\beta}_M = (X_M' X_M)^{-1} X_M'y$.
The working model $M$ is correct
if $X_M \beta_M = \mu$; in that case, we have $v_M = u$.
Otherwise, that is, if $X_M \beta_M \neq\mu$, the working model is incorrect,
and we have $v_M = \mu- X_M \beta_M + u$.
Irrespective of whether the working model is correct or not,
we always have $\hat{\beta}_M \sim N(\beta_M, \sigma^2
(X_M'X_M)^{-1} )$;
in particular, $\hat{\beta}_M$ is an unbiased estimator for $\beta_M$,
irrespective of whether or not the model $M$ is correct.
As noted earlier, we assume that the variance estimator $\hat{\sigma}^2$
is independent of the collection of estimators
$\hat{\beta}_M$ for $M \in\mathcal{M}$.

To pinpoint the regression coefficient of a given regressor $X_j$
in a model $M$ it appears in, we write
$\beta_{j\cdot M}$ for that component of $\beta_{M}$
that corresponds to the regressor $X_j$ for each $j \in M$.
Similarly, the components of $\hat{\beta}_M$ are indexed
as $\hat{\beta}_{j\cdot M}$ for $j\in M$.
This convention is called ``full model indexing'' in \citet{Ber13a}.

Consider now a model selection procedure, that is, a data-driven rule
that selects a model $\hat{M} \in\mathcal{M}$ from the pool $\mathcal{M}$ of
candidate models
and the resulting post-model-selection estimator $\hat{\beta}_{\hat{M}}$.
The\vspace*{1pt} coverage target considered in \citet{Ber13a}
is $\beta_{\hat{M}}$ or components thereof.
Note that this coverage target is random,
because it depends on the outcome of the model selection procedure.

\begin{remark}\label{r2}
(i) At least one author of the present paper
believes that the merits
of $\beta_{\hat{M}}$
as a coverage target for inference are debatable: For example,
the meaning of the first coefficient of $\beta_{\hat{M}}$ depends on
the selected model and hence also on the training data $(y,X)$;
the same
applies to the dimension of $\beta_{\hat{M}}$.
In particular, we stress that different model selection
procedures (e.g., AIC, BIC, the LASSO, etc.) lead to different
targets $\beta_{\hat{M}}$.
We refer to \citet{Ber13a} for further discussion and motivation
for studying $\beta_{\hat{M}}$.
These authors make the case for $\beta_{\hat{M}}$ by arguing that
the relevant setting is one where no correct overall model is
available; however, in this situation the subsequent remark becomes
especially important.

(ii) While the model \eqref{model} is nonparametric,
the distributional requirements on $\hat{\sigma}^2$
obviously are rather
restrictive.
However, these are the assumptions underlying the analysis in
\citet{Ber13a}, and we adopt them here in order to be in line with
that reference. A leading case where
these requirements are fulfilled is when \eqref{model}
is replaced by the \textit{parametric} model $y = X \beta+ u$,
when $X$ is as before and is assumed to be of full column rank $p < n$,
and when $\hat{\sigma}^2$ is the usual unbiased
variance estimator in that model and $r$ is set to $n-p$.
In this leading case, however,
the true parameter $\beta$ in the overall model
is well-defined and will then typically be the prime target of
statistical inference, rather than the nonstandard coverage target
introduced in \citet{Ber13a}.
Outside of the parametric model just discussed,
the requirements on $\hat{\sigma}^2$ made in \citet{Ber13a},
and also here, will only be satisfied in certain special cases, some
of which are discussed at the end of Section~2.2 in \citet{Ber13a}.
[The requirements on $\hat{\sigma}^2$ are also fulfilled
(with $r=n-q$), if we would maintain a true parametric model
$y=Z \theta+ u$ for some observed $n\times q$ matrix $Z$
of rank $q<n$ that contains $X$ as a submatrix;
however, in this case
one is back to the leading case discussed above, after
redefining $\mathcal{M}$ appropriately.]
\end{remark}

In this paper, we will mainly focus on confidence intervals for the
coefficient of one particular regressor in the selected model. Without
loss of generality, assume that $X_1$ is the regressor of interest
and that the coverage target is $\beta_{1\cdot\hat{M}}$. To ensure
that this quantity is always well-defined, we assume that the first
regressor $X_1$ is contained in all candidate models under consideration,
that is, we assume that $1 \in M$ for each $M \in\mathcal{M}$.
We seek to construct
confidence intervals for $\beta_{1\cdot\hat{M}}$ that are
of the form
\[
\hat{\beta}_{1\cdot\hat{M}} \pm K \hat{\sigma}_{1\cdot\hat{M}}
\]
for some constant $K>0$,
with $\hat{\sigma}^2_{1\cdot M}$ defined by
$\hat{\sigma}^2_{1\cdot M} = \hat{\sigma}^2 [(X_M'X_M)^{-1}]_{1,1}$,
where $[\cdots]_{1,1}$ denotes the first diagonal element of the indicated
matrix.
Here, we abuse notation and write $a\pm b$ for the interval
$[a-b, a+b]$.
For a given level $1-\alpha$ with $0<\alpha<1$,
the constant $K$ should be chosen such
that the minimal coverage probability is at least $1-\alpha$, that is,
such that
%
\begin{equation}
\label{CI}
\hspace*{9pt}\inf_{\mu,\sigma} \mathbb{P}_{\mu,\sigma} (
\beta_{1\cdot\hat{M}} \in \hat {\beta}_{1\cdot\hat{M}} \pm K\hat{
\sigma}_{1\cdot\hat{M}} ) \geq 1-\alpha.
\end{equation}

Because the distribution of
$(\hat{\beta}_{1\cdot M} -\beta_{1\cdot M})/ \hat{\sigma}_{1\cdot M}$
is independent of unknown parameters and also independent of $M$,
it follows, for \textit{fixed} $M$, that a confidence interval for
$\beta_{1\cdot M}$ with minimal coverage probability $1-\alpha$
is given by the textbook interval
$\hat{\beta}_{1\cdot M} \pm K_N \hat{\sigma}_{1\cdot M}$,
where $K_N$ is the $(1-\alpha/2)$-quantile of the distribution of
$(\hat{\beta}_{1\cdot M}-\beta_{1\cdot M}) / \hat{\sigma}_{1\cdot
M}$---a standard normal
distribution in the known-variance case and a $t$-distribution with
$r$ degrees of freedom in the unknown-variance case.
In view of this,
it is tempting to consider, as a confidence interval for
$\beta_{1\cdot\hat{M}}$, the interval
$\hat{\beta}_{1\cdot\hat{M}} \pm K_N \hat{\sigma}_{1\cdot\hat{M}}$.
Because this construction ignores the model selection step and treats
the selected model $\hat{M}$ as fixed, we will call this the ``naive''
confidence interval.

The PoSI-interval developed in \citet{Ber13a} is obtained
by first constructing simultaneous
confidence intervals for the components of $\beta_M$
that are centered at the corresponding components of $\hat{\beta}_M$,
for each $M \in\mathcal{M}$,
with coverage probability $1-\alpha$:
More formally, the PoSI-constant $K_P$ is the unique solution to
%
\begin{eqnarray}
&& \inf_{\mu,\sigma} \mathbb{P}_{\mu,\sigma} (
\beta_{j\cdot M} \in\hat{\beta }_{j\cdot M} \pm K_P \hat{
\sigma}_{j\cdot M} \dvtx j \in M, M \in\mathcal{M} )\hspace*{-10pt} \nonumber
\\[-3pt]
\label{KP}
\\[-8pt]
\nonumber
&&\quad= 1-\alpha,\hspace*{-10pt}
\end{eqnarray}
where the quantities $\hat{\sigma}_{j\cdot M}^2$ are defined like
$\hat{\sigma}_{1\cdot M}^2$ but with $j$ replacing $1$.
By construction, the PoSI-constant $K_P$ is such that
we obtain simultaneous confidence intervals for the components
of $\beta_{\hat{M}}$ that are centered at the corresponding components
of $\hat{\beta}_{\hat{M}}$. In other words, \eqref{KP} implies
%
\begin{eqnarray}
&& \inf_{\mu,\sigma} \mathbb{P}_{\mu,\sigma} (
\beta_{j\cdot\hat{M}} \in\hat {\beta}_{j\cdot\hat{M}} \pm K_P \hat{
\sigma}_{j\cdot\hat{M}} \dvtx j \in\hat{M} )
\nonumber
\\[-8pt]
\label{PoSI}
\\[-8pt]
\nonumber
&&\quad \geq 1-\alpha.
\end{eqnarray}
In particular, \eqref{CI} holds when $K_P$ replaces $K$.
For computing the constant $K_P$, we note that
the probability in \eqref{KP} can also be written as
$\mathbb{P}_{\mu,\sigma}( | \hat{\beta}_{j\cdot M} - \beta
_{j\cdot M}| /
\hat{\sigma}_{j\cdot M} \leq K_P\dvtx
j \in M, M \in\mathcal{M})$.
This probability is not hard to compute,\vspace*{1pt}
because it involves only the random variables
$(\hat{\beta}_{j\cdot M} - \beta_{j\cdot M}) / \hat{\sigma
}_{j\cdot M}$,
which are
(dependent) standard normal in the known-variance case
and (dependent) $t$-distributed in the unknown-variance case,
with an obvious dependence structure only depending on $X$.
In particular, the probability in \eqref{KP} does not
depend on $\mu$ or $\sigma^2$. Similar considerations apply, mutatis
mutandis, to the constant $K_{P1}$ that is introduced in the following
paragraph.

A modification of the preceding
procedure, which is also proposed in~\citet{Ber13a},
is useful when inference is focused on a particular component of
$\beta_{\hat{M}}$, instead of on all components.
Recall that the coverage target in \eqref{CI} is the first component
of $\beta_{\hat{M}}$, that is, $\beta_{1\cdot\hat{M}}$.
The PoSI1-constant $K_{P1}$ provides simultaneous confidence intervals
for $\beta_{1\cdot M}$ centered at $\hat{\beta}_{1\cdot M}$ for each
$M\in\mathcal{M}$.
In particular, $K_{P1}$ is the unique solution to
%
\begin{eqnarray}
&& \quad\inf_{\mu,\sigma} \mathbb{P}_{\mu,\sigma} (
\beta_{1\cdot M} \in\hat{\beta }_{1\cdot M} \pm K_{P1} \hat{
\sigma}_{1\cdot M} \dvtx M\in\mathcal{M} )
\nonumber
\\[-8pt]
\label{PoSI1}
\\[-8pt]
\nonumber
&&\quad \quad= 1-\alpha.
\end{eqnarray}
Again, by construction, \eqref{CI} holds when $K_{P1}$ replaces~$K$.

Like the PoSI-constants discussed so far, other procedures
for controlling the family-wise error rate can be used.
Consider, for example,
Scheff\'{e}'s method: Recall that $X$ denotes the matrix of all available
explanatory variables, and note that
$(\hat{\beta}_{j\cdot M} - \beta_{j\cdot M})$ is a linear function
of $Y-\mu$,
that is, a function of the form $\upsilon' (Y-\mu)$, for a certain vector
$\upsilon\neq0$
in the span of $X$. The Scheff\'{e} constant $K_{S}$ is chosen such that
\[
\mathbb{P}_{\mu,\sigma} \biggl( \mathop{\sup_{\nu\in\operatorname{span}(X)}}_{\nu\neq0} \frac{ \nu' (Y-\mu)}{\hat{\sigma}\|\nu\|}
\leq K_{S} \biggr) = 1-\alpha.
\]
Then the relations \eqref{PoSI} and, in particular,
\eqref{CI} hold when $K_S$ replaces both $K$ and $K_P$.
Note that the probability in the preceding display does not depend
on $\mu$ and $\sigma$, and that the constant $K_S$ is easily computed as
follows:
Let $s$ denote the rank of $X$.
In the known-variance case, $K_S$ is the square root of
the $(1-\alpha)$-quantile of a chi-square distribution with $s$ degrees
of freedom. In the unknown-variance case, $K_S$ is the square root of
the product of $s$ and the $(1-\alpha)$-quantile of an $F$-distribution
with $s$ and $r$ degrees of freedom.

Using the constants $K_P$, $K_{P1}$, or $K_S$ gives valid confidence intervals
post-model-selection, that is, intervals that satisfy \eqref{CI},
because these constants
give simultaneous confidence intervals
for all quantities of interest that can occur; for example,
\eqref{PoSI} follows from \eqref{KP}, which in turn
guarantees that \eqref{CI} holds when $K_P$ replaces $K$.
One advantage of this is that a coverage
probability of at least $1-\alpha$
is guaranteed,
\textit{irrespective} of the model selection procedure $\hat{M}$
(as long as it takes values in $\mathcal{M}$).
In particular, this is guaranteed even if the model is selected
by statistically inane methods like the SPAR-procedure mentioned
in Section~4.9 of \citet{Ber13a}.
The price for this is that the PoSI constants $K_P$ and $K_{P1}$
may be overly conservative for a \textit{particular}
model selection procedure $\hat{M}$.
[In this context, we note that equality holds in \eqref{PoSI}
for the SPAR-procedure,
and that equality holds in \eqref{CI}
for a variant of the SPAR-procedure which selects that model
$\hat{M}$ which maximizes $|\hat{\beta}_{1\cdot M}|/\hat{\sigma
}_{1\cdot M}$
over $M\in\mathcal{M}$.
Because such model selection procedures are hard to justify from a
statistical perspective, we will not further consider SPAR and its variant
here.]

Last, we will also consider the obvious approach where one
chooses the smallest constant $K$ such that \eqref{CI} is satisfied.
We will denote this
constant by $K_\ast$ (provided it exists).
This is, of course, a well-known standard construction; see
\citet{Bic77a}, page~170, for example.
By definition, the interval in \eqref{CI} with $K_\ast$ replacing $K$
is the shortest interval of that form
whose minimal coverage probability is $1-\alpha$.
Note that $K_\ast$ depends on the model selection procedure in
question, and that computation of this quantity
can be cumbersome as it requires computation of the
finite-sample distribution of
$\hat{\beta}_{1\cdot\hat{M}} / \hat{\sigma}_{1\cdot\hat{M}}$.
However, explicit computation of this constant is feasible in some
cases [cf. the results in Section~\ref{theory} and also
the more general results of \citet{Lee02d}],
and this constant can also be computed or approximated
in a variety of other scenarios [e.g., by
adapting the results of \citet{Poe10a} or the
procedures of \citet{And09a}].
Also note that we have $K_\ast\leq K_{P1} \leq K_P \leq K_S$
by construction.

The procedures discussed so far are concerned with
coverage targets like $\beta_{\hat{M}}$ that depend on the
selected model.
This should be compared to the more classical parametric setting where
the coverage target is the underlying true parameter:
Assume that the data is generated by an overall linear model, that is,
assume that
the parameter $\mu$ in \eqref{model} satisfies
$\mu= X \beta$ for the overall regressor matrix $X$ introduced earlier,
and that $\operatorname{rank}(X) = p<n$ holds.
And assume that inference is focused on (components of) the parameter
$\beta$.
In this setting, the effect of model selection on subsequently
constructed confidence intervals can be dramatic. For example,
\citet{Kab06a} show that the minimal coverage probability
of the ``naive'' confidence interval for $\beta_1$, that is, the quantity
%
\[
\inf_{\beta,\sigma} \mathbb{P}_{X\beta,\sigma} (
\beta_1 \in \hat{\beta }_{1\cdot\hat{M}} \pm K_N \hat{
\sigma}_{1\cdot\hat{M}} ),
\]
can be much smaller than the nominal coverage probability $1-\alpha$;
in fact, this minimal coverage probability
can, for example, be smaller than 0.5, depending on the regressor
matrix $X$
in the overall model $y = X\beta+ u$.
The main reason for this more dramatic effect is that
$\hat{\beta}_{1\cdot M}$ is a biased estimator for $\beta_1$
whenever the
model $M$ is incorrect, whereas $\hat{\beta}_{1\cdot M}$ is always unbiased
for $\beta_{1\cdot M}$.
Of course, valid confidence intervals post-model-selection can also
be constructed when the coverage target is $\beta_1$, namely, by
replacing $K_N$
in the preceding display by the smallest constant $K$ such that
the resulting minimal coverage probability equals $1-\alpha$
(provided it exists).
For the computation or approximations of this constant in
particular situations, we refer to the papers cited in the preceding
paragraph.

\section{Some Finite-Sample Results}
\label{theory}

In this section we give a finite-sample analysis of the
confidence intervals discussed so far, where we consider
a simple model selection
procedure that selects among two nested models using a likelihood ratio test.
More precisely, maintaining the setting of Section~\ref{setting},
let $X$ now be an $n\times2$ matrix of rank 2, and
assume that $\mathcal{M} = \{M_1, M_2\}$ with $M_1 = \{1\}$ and
$M_2 = \{1,2\}$ throughout this section.
For the model selector, we
set $\hat{M} = M_2$ if $|\hat{\beta}_{2\cdot M_2} |/ \hat{\sigma
}_{2\cdot M_2}$
is larger than~$C$, and $\hat{M}=M_1$ otherwise, where $C>0$ is a
user-specified constant.
Arguably, any reasonable model selection
procedure in this setting must be equivalent to a likelihood ratio
test, at least asymptotically; cf. \citet{Kab06a}.
In the numerical examples that follow,
we will consider $C=\sqrt{2}$, such that the resulting model
selector $\hat{M}$ corresponds to selection by
the classical Akaike information criterion (AIC);
this model selector is asymptotically equivalent to
several other model selectors,
including the GCV model selection criterion of \citet{Cra78a}
and the $S_p$ criterion of \citet{Tuc67a}; cf. \citet{Lee05a}.
Furthermore, we will also consider
$C=\sqrt{\log(n)}$, corresponding to the BIC model selection criterion.
Throughout this section, let $\phi(\cdot)$ and $\Phi(\cdot)$
denote the density and the cumulative distribution function (c.d.f.)
of the univariate standard Gaussian distribution,
and set $\Delta(x,c) = \Phi(x+c) - \Phi(x-c)$.
And, last, we will write $\rho$ for
the correlation coefficient between
the two components of $\hat{\beta}_{M_2}$, that is,
$\rho= -[(X_{M_2}'X_{M_2})^{-1}]_{1,2} (
[(X_{M_2}'X_{M_2})^{-1}]_{1,1}
[(X_{M_2}' \cdot\break X_{M_2})^{-1}]_{2,2}
)^{-1/2}$.

The following result describes the
coverage probability of the interval
$\hat{\beta}_{1\cdot\hat{M}} \pm K \hat{\sigma}_{1\cdot\hat{M}}$
in two scenarios, namely, when the coverage target is
$\beta_{1\cdot\hat{M}}$ and when the coverage target is
$\beta_{1\cdot M_2}$.
Note that, in case the model $M_2$ is correct, that is, if we have $\mu
= X \beta$
for some $\beta\in\mathbb{R}^2$,
and hence also $y = X \beta+ u$, then this second scenario
reduces to the classical parametric setting described at
the end of Section~\ref{setting}; in particular, we then have
$\beta_{M_2} = \beta$
and, thus, $\beta_{1\cdot M_2} = \beta_1$.

\begin{proposition}
\label{P1}
In the setting of this section, we have
\begin{eqnarray*}
&&\mathbb{P}_{\mu,\sigma} ( \beta_{1\cdot\hat{M}} \in \hat{
\beta}_{1\cdot\hat{M}} \pm K \hat{ \sigma}_{1\cdot\hat{M}} )
\\
&&\quad=
\mathbb{E} \biggl[ \Delta \biggl(0, \frac{\hat{\sigma}}{\sigma}K \biggr) \Delta \biggl(
\zeta,\frac{\hat{\sigma}}{\sigma}C \biggr)\\
&&\qquad \hspace*{12pt}{}+ \int_{-({\hat{\sigma}}/{\sigma})K}^{({\hat{\sigma
}}/{\sigma})K}
\biggl(1- \Delta \biggl( \frac{\zeta+ \rho z}{\sqrt{1-\rho^2}}, \frac{({\hat{\sigma}}/{\sigma})C}{\sqrt{1-\rho^2}} \biggr) \biggr)
\\
&&\hspace*{120pt}\qquad\qquad\qquad{}\cdot\phi(z) \,d z \biggr]
\end{eqnarray*}
and
\begin{eqnarray*}
&& \mathbb{P}_{\mu,\sigma} ( \beta_{1\cdot M_2} \in \hat{
\beta}_{1\cdot\hat{M}} \pm K \hat{ \sigma}_{1\cdot\hat{M}} ) \\
&&\quad= \mathbb{P}_{\mu,\sigma}
( \beta_{1\cdot\hat{M}} \in \hat{\beta}_{1\cdot\hat{M}} \pm K \hat{
\sigma}_{1\cdot\hat{M}} )
\\
&& \qquad{}+\mathbb{E} \biggl[ \biggl( \Delta \biggl(\frac{\rho\zeta}{\sqrt{1-\rho^2}},
\frac{\hat{\sigma}}{\sigma}K \biggr) \\
&&\qquad\qquad\hspace*{32pt}{}- \Delta \biggl(0, \frac{\hat{\sigma}}{\sigma}K \biggr) \biggr)
\Delta \biggl( \zeta,\frac{\hat{\sigma}}{\sigma}C \biggr) \biggr],
\end{eqnarray*}
with $\zeta= \beta_{2\cdot M_2}/\operatorname{SD}(\hat{\beta}_{2\cdot M_2})$,
where $\operatorname{SD}(\cdot)$ denotes the standard deviation.
The expectations on the right-hand sides are taken with respect to
$\hat{\sigma}/\sigma$. In the known-variance case, $\hat{\sigma
}/\sigma$
is constant equal to one and the expectations are trivial;
in the unknown-variance case, $\hat{\sigma}/\sigma$
is distributed like the square root of a chi-squared distributed random
variable with $r$ degrees of freedom divided by $r$, that is,
$\hat{\sigma}/\sigma\sim\sqrt{\chi^2_r/r}$.
\end{proposition}

\begin{pf}
The statements for the known-variance case are simple adaptations of the
finite-sample statements of Proposition~3 in \citet{Kab06a}.
For the unknown-variance case, it suffices to note that $\hat{\sigma
}/\sigma$
is independent of $\{\hat{\beta}_{M_1}, \hat{\beta}_{M_2}\}$. With this,
the statements are then obtained by conditioning on $\hat{\sigma
}/\sigma$
and by using the formulae for the known-variance case derived earlier.
\end{pf}

Proposition~\ref{P1} provides explicit formulas that
also allow us to compute (minimal) coverage probabilities numerically.
For the following discussion, fix the values of
$C$ and $K$, that is, the critical value $C$
of the hypothesis test that is used for model selection and
the value $K$ that governs the length of the confidence interval
post-model-selection.
We first\vspace*{1.5pt} note that
$\mathbb{P}_{\mu,\sigma}(\beta_{1\cdot M_2} \in\hat{\beta
}_{1\cdot\hat{M}}\pm
K \hat{\sigma}_{1\cdot\hat{M}})$
is strictly smaller than
$\mathbb{P}_{\mu,\sigma}(\beta_{1\cdot\hat{M}} \in\hat{\beta
}_{1\cdot\hat{M}}\pm
K \hat{\sigma}_{1\cdot\hat{M}})$ whenever $\rho\zeta\neq0$, because
the two probabilities differ by a correction term (namely, the expected
value on the right-hand side of the second display in Proposition~\ref{P1})
which is negative whenever $\rho\zeta\neq0$.
If $\rho\zeta= 0$, the two probabilities are equal.
And if $\rho=0$, it is easy to see that both probabilities are equal to
$\mathbb{E}[\Delta(0,K\hat{\sigma}/\sigma)] =
F(K)-F(-K)$, irrespective of $\zeta$,
where $F$ denotes the c.d.f. of a $t$-distribution with $r$ degrees of
freedom in the unknown-variance case and the standard Gaussian c.d.f.
in the known-variance case.
Next, we note that the coverage probabilities depend only on
$r$, $\zeta$, and $\rho$.
(Recall that $r$ denotes the degrees of freedom of $\hat{\sigma}^2$ in
the unknown-variance case, and that we have set $r=\infty$ in the
known-variance case.)
Note that $\zeta$ is a function of the regressor matrix $X_{M_2}$
and of the unknown parameters $\mu$ and $\sigma^2$, while
$\rho$ is a function of $X_{M_2}$ only.
Moreover, it is easy to
see that the coverage probabilities are symmetric both in $\zeta$
and in $\rho$ around the origin.
Concerning the influence of $r$,
it can be shown that the coverage probabilities for the known-variance
case provide a uniform approximation to those in the
unknown-variance case, uniformly in the unknown parameters,
where the approximation error goes to zero as $r\to\infty$;
this follows from the results of \citet{Lee02d} using standard arguments.
In the examples that follow, we found that the results for the
known-variance case and for the unknown-variance case are similar,
and that these results are visually hard to distinguish from each other,
unless $r$ is extremely small like, for example, 3.
We therefore focus on the known-variance case in the following because
it provides a good approximation to the unknown-variance case
as long as $r$ is not too small.

\begin{figure}

\includegraphics{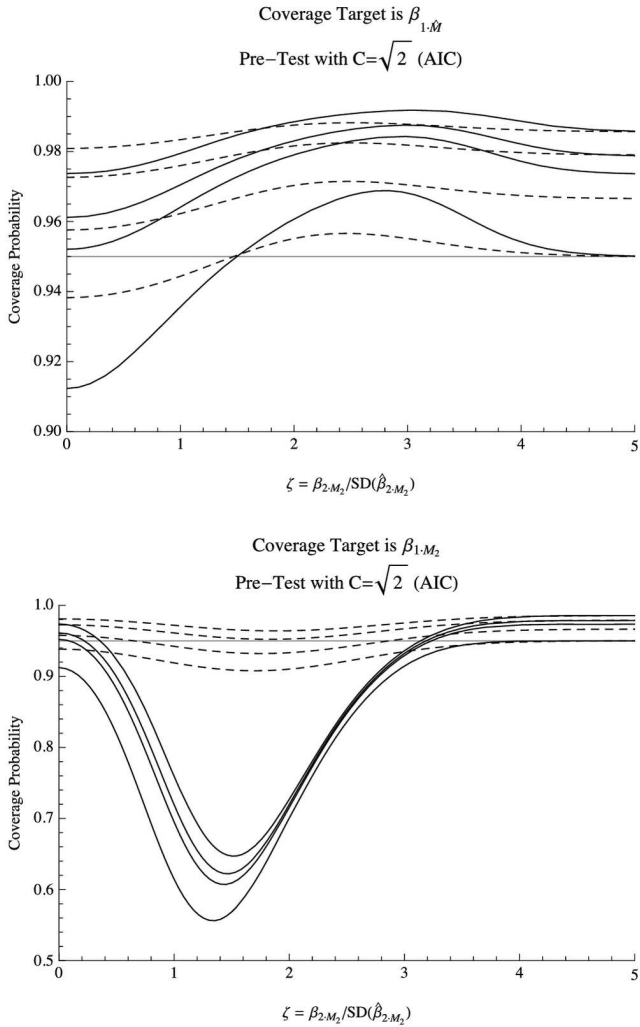}

\caption{Coverage probability of several confidence intervals in
the known-variance case, as a function
of the scaled parameter
$\zeta= \beta_{2\cdot M_2} / \operatorname{SD}(\hat{\beta}_{2\cdot M_2})$,
using the model selection procedure with $C=\sqrt{2}$, that is, AIC.
The nominal coverage probability is $1-\alpha= 0.95$, indicated
by a gray horizontal line.
The coverage target is $\beta_{1\cdot\hat{M}}$ (top panel)
and $\beta_{1\cdot M_2}$ (bottom panel). In each panel,
the four solid curves are computed for $\rho=0.9$,
and the four dashed curves are for $\rho=0.5$. The curves in each
group of four are ordered: Starting from the top, the curves show
the coverage probabilities
for $K_S$ (Scheff\'{e}), $K_P$ (PoSI), $K_{P1}$ (PoSI1), and $K_N$ (naive).}\label{fig1}
\end{figure}

We proceed to comparing the case where the coverage target is
$\beta_{1\cdot\hat{M}}$ as in \citet{Ber13a} with the more standard
case where the
coverage target is the parameter $\beta_{1\cdot M_2}$, in terms of
the coverage probabilities of confidence intervals post-model-selection.
Recall that the nonstandard target depends on the training data
as well as on the model selection procedure employed, whereas the
standard target does not.
Consider first the case where $C=\sqrt{2}$, corresponding to the
AIC model selector.
For several of the confidence intervals introduced in the preceding section,
the results are visualized in Figure~\ref{fig1}, for the case where the
coverage target is $\beta_{1\cdot\hat{M}}$ (top panel) and for the
case where the coverage target is $\beta_{1\cdot M_2}$
(bottom panel). Note that the range of the vertical axes
(displaying coverage probability) in the two panels is quite different.

In each panel of Figure~\ref{fig1},
we see that the effect of model selection on the resulting
coverage probabilities depends on the correlation coefficient $\rho$,
with larger values of $\rho$ corresponding to
smaller minimal coverage probabilities. But the strength of the effect
varies greatly with the scenario, that is, on whether the coverage target
is $\beta_{1\cdot\hat{M}}$ or $\beta_{1\cdot M_2}$.
When the coverage target is
$\beta_{1\cdot\hat{M}}$ (top panel in Figure~\ref{fig1}), we see that the
effect of
model selection is comparatively minor:
The smallest coverage probabilities are always obtained for the
``naive'' interval, whose coverage probability here can be
smaller as well as larger than the nominal 0.95. Irrespective of the
true parameters, the actual coverage
probability of the ``naive'' interval is quite close to the
nominal one here. The other intervals, that is, the PoSI1-, the PoSI-, and
the Scheff\'{e}-interval, all have coverage probabilities larger than $0.95$.
(The minimal coverage probabilities here are obtained for $\zeta=0$,
but we found this not to be the case for other model selection procedures,
i.e., for other values of~$C$.)
When the coverage target is $\beta_{1\cdot M_2}$
(bottom panel in Figure~\ref{fig1}), however, we get a very different picture:
For $\rho=0.9$, the minimal coverage probability of all the
intervals considered there is much smaller than 0.95, with minima
between 0.55 (``naive'') and 0.65 (Scheff\'{e}).
For $\rho=0.5$, the minimal coverage probabilities of the ``naive'' interval
and of the PoSI1-interval are below, while those of the other intervals
are above, the nominal 0.95.
For very small values of $\rho$, the coverage probabilities of
all the intervals considered in Figure~\ref{fig1} are visually indistinguishable
from horizontal lines as a function of $\zeta$ (and hence are not
shown here),
irrespective of the coverage target.
For $\rho=0.1$, for example,
the coverage probability of
the ``naive'' interval is about 0.95, while that of the
other intervals is above $0.95$, ordered by their length.
(This should not come as a surprise since in case $\rho= 0$
model selection has no effect on estimating the regression coefficients;
furthermore, the two targets are identical in this case.)

Figure~\ref{fig1} illustrates that the coverage probability of
confidence intervals post-model-selection depends crucially
on whether the coverage target is $\beta_{1\cdot\hat{M}}$ as
in \citet{Ber13a} or the more classical coverage target $\beta_{1\cdot M_2}$.
We stress here again that the PoSI-intervals and the Scheff\'{e}-interval
have not been designed to deal with the case where the coverage target
is $\beta_{1\cdot M_2}$.
For a more detailed analysis of the ``naive'' interval in
the case where the coverage target
is $\beta_{1\cdot M_2}$, we refer to \citet{Kab06a}.

For the other values of $C$ that we consider, that is, for
$C=\sqrt{\log(n)}$ for various values of $n$, we found the following:
When the coverage target is $\beta_{1\cdot\hat{M}}$, the results
are very similar to those shown in the top panel of Figure~\ref{fig1}.
To conserve space, we do not show these results here.
When the target is $\beta_{1\cdot M_2}$, the resulting curves
are of the same shape but steeper, with coverage probabilities
decreasing as $C$ increases.
This is so because larger values of $C$ lead to more frequent selection
of the smaller model $M_1$, causing more bias in the resulting
post-model-selection estimator; we refer
to \citet{Lee03a} and, in particular,
Figure~3 in that reference, for further discussion and analysis of this
phenomenon.

We next compare the confidence intervals for $\beta_{1\cdot\hat{M}}$
introduced in Section~\ref{setting} through their minimal coverage
probability as a function of the correlation coefficient $\rho$.
In particular, for various values of $C$,
we compute the quantity on the left-hand side of \eqref{CI}
for specific $K$'s, namely, for $K_N$ (``naive''),
for $K_P$ (PoSI), for $K_{P1}$ (PoSI1), for $K_S$ (Scheff\'{e}),
and for $K_\ast$ (the smallest valid $K$).
By construction, we have $K_\ast\leq K_{P1} \leq K_P \leq K_S$,
so that the resulting curves of minimal coverage probabilities
are also arranged in increasing order.
\begin{figure}[b]

\includegraphics{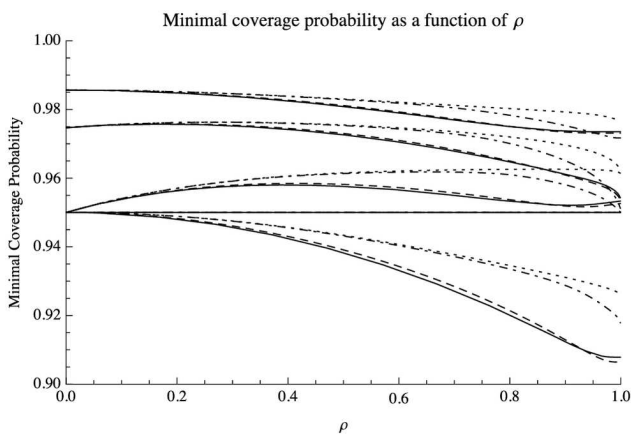}

\caption{Minimal coverage probabilities of the confidence intervals
for $\beta_{1\cdot\hat{M}}$ as a function of
$\rho$ in the known-variance case,
for $C=\sqrt{2}$ (solid curves), $C=\sqrt{\log(10)}$ (dashed curves),
$C=\sqrt{\log(100)}$ (dot-dashed curves),
and $C=\sqrt{\log(1000)}$ (dotted curves).
The nominal coverage probability is $1-\alpha= 0.95$.
For each value of $C$, the corresponding
five curves are ordered: Starting from the top,
the curves correspond to the intervals with $K_S$, $K_P$,
$K_{P1}$, $K_\ast$, and $K_N$.}\label{fig2}
\end{figure}

All the minimal coverage probabilities shown in Figure~\ref{fig2} are within
5\% of the nominal level 0.95.
For the ``naive'' intervals corresponding to $K_N$ (the first four curves
from the bottom), the minimal
coverage probability is below 0.95 (except for the trivial case where
$\rho=0$), but not by much.
The intervals with $K_\ast$ have minimal coverage probabilities
of exactly 0.95, for every value of $C$, by construction
(but note that $K_\ast$ depends on $C$, whereas $K_S$, $K_P$, $K_{P1}$,
and $K_N$ do not).
Hence, the curves corresponding to the $K_\ast$'s for the four values
of $C$
considered here are constant and sit on top of each other.
And, again by construction, all other
intervals are slightly too large in the sense that their minimal coverage
probability exceeds the nominal level 0.95.
Concerning the influence of $C$, we see that larger values of $C$
correspond to slightly larger minimal coverage probabilities for
the intervals corresponding to $K_N$, $K_{P1}$, $K_P$, and $K_S$,
and for most values of $\rho$;
it should be noted, however, that---in contrast to the case of the
standard target---here the target changes with $C$.
Overall, the difference between the coverage probabilities of
all these intervals is not dramatic.

Last, we compare the confidence intervals for $\beta_{1\cdot\hat{M}}$
through the values of the constants $K$ that correspond
to the intervals in question.
By construction, $K_S$ and $K_N$ are constant as a function of $\rho$.
Note that the constants
$K_N$, $K_P$, $K_{P1}$, and $K_S$ do not depend on the model
selection procedure that is being used (and thus not on $C$),
while the constant $K_\ast$ does
depend on the model selection procedure (and thus on $C$).
For a given model selection procedure,
the constant $K_\ast$ is the smallest number $K$
for which \eqref{CI} holds; in particular,
the interval corresponding to $K$ has minimal coverage
probability smaller/equal/larger than $1-\alpha$ if and only if
$K$ is smaller/equal/larger than $K_\ast$.

\begin{figure}

\includegraphics{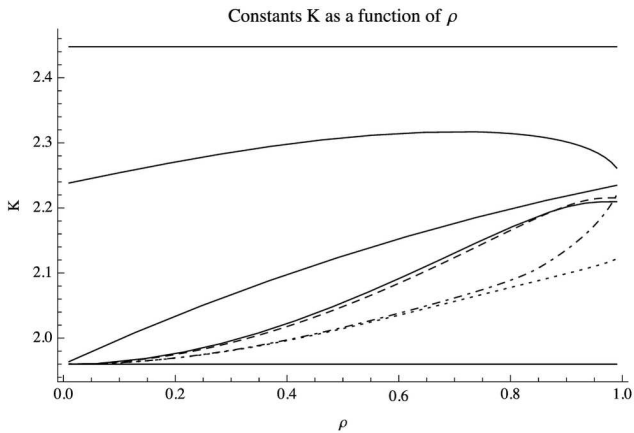}

\caption{The constants $K$ that govern the width of the confidence intervals
as a function of $\rho$ in the known-variance case,
using the model selection procedure with critical value $C$.
The nominal coverage probability is $1-\alpha= 0.95$.
Starting from the top, the five solid
curves show $K_S$, $K_P$, $K_{P1}$, $K_\ast$
for $C=\sqrt{2}$ (AIC), and $K_N$.
The remaining curves show $K_\ast$ for $C=\sqrt{\log(10)}$ (dashed curve),
for $C=\sqrt{\log(100)}$ (dash-dotted curve), and
for $C=\sqrt{\log(1000)}$ (dotted curve).}\label{fig3}
\end{figure}

The interpretation of Figure~\ref{fig3} is similar to that of Figure~\ref{fig2},
the main difference being that the lengths considered here
are somewhat more distorted than the minimal coverage probabilities
considered earlier. The ``naive'' interval is up to about 10\% too short,
while the intervals corresponding to $K_{P1}$, $K_P$, and $K_S$ are
too long, namely, by up to about
5\%, 15\%, 25\%, respectively.
We also see that $K_\ast$ decreases as $C$ increases for most values
of $\rho$, which is consistent with the observations made in
the second-to-last paragraph.

\section{Simulation Study}
\label{simulation}

We now compare the ``naive'' interval, the PoSI1 interval, and
(a variant of) the PoSI interval for $\beta_{1\cdot\hat{M}}$
by their respective minimal coverage probabilities
in a simulation study
where the data are generated from a Gaussian overall linear model
$M_{\mathrm{full}}$, say, of the form
$Y=X\beta+ u$
with 30 observations, 10 explanatory variables, and i.i.d. standard normal
errors.
Moreover, we also study these intervals when the coverage target
is $\beta_1 = \beta_{1\cdot M_{\mathrm{full}}}$
(instead of $\beta_{1\cdot\hat{M}}$).
For the estimator~$\hat{\sigma}^2$, we use the usual unbiased variance
estimator obtained by fitting the overall model; hence, we have $r=n-p=20$
here.
[To be precise, while the constants $K_N$ as well as $K_{P1}$ are
computed as detailed in Section~\ref{setting}, we consider
instead of $K_P$ defined by \eqref{KP} the larger constant
$K_{P'}$ which is obtained from \eqref{KP} when $\mathcal{M}$
is replaced by the collection of \textit{all} nonempty subsets
of $\{1,\ldots, p\}$. We shall refer to the resulting interval
also as a PoSI-interval in this section.
The reason for this choice is that code for computing
$K_{P'}$ is publicly available from the authors of \citet{Ber13a},
so that $K_{P'}$ is the PoSI-constant
likely to be used by practitioners.
Note that $K_{P1} \leq K_P \leq K_{P'}$ holds, and hence the
performance of the interval based on $K_P$ can be easily deduced
from Table~\ref{tab1}.]

As model selectors, we consider AIC, BIC, and
the LASSO: For AIC we use the \texttt{step()} function in
R with its default settings, subject to the constraint that the
regressor of interest, that is, the first one, is always included;
this corresponds to minimizing the
AIC objective function through a greedy general-to-specific search
over the $2^{9}$ candidate models (i.e., $\mathcal{M}$ consists
of all submodels of the overall model that contain the
first regressor).
Similarly, for BIC we use the \texttt{step()} function
with the penalty parameter equal to $\log(30)$.
And for the LASSO, we basically select those regressors for which
the LASSO-estimator has nonzero coefficients. [More precisely,
we use the \texttt{lars} package in R and
follow suggestions outlined in \citet{Efr04a}, Section 3.4:
To protect the regressor of interest (the first one), we
first compute the residual of the orthogonal projection of $y$ on the
first regressor; write $\tilde{y}$ for this residual vector, and
write $\tilde{X}$ for the regressor matrix $X$ with the first column removed.
We then compute\vspace*{1pt} the LASSO-estimator for a regression of $\tilde{y}$ on
$\tilde{X}$ using the \texttt{lars()} function;
the LASSO-penalty is chosen by 10-fold cross-validation
using the \texttt{cv.lars()} function (in both functions, we set the
\texttt{intercept} parameter to \texttt{FALSE}, and otherwise use
the default settings).
The selected model is comprised of those regressors in $\tilde{X}$ for which
the corresponding LASSO coefficients are nonzero, plus the first column
of $X$.]

Three designs are considered
for the design matrix $X$:
For design 1, we take the regressor matrix from
the data-example from Section~3 of \citet{Kab06a} (for which the minimal
coverage probability of a ``naive'' nominal 95\% interval for $\beta_1$,
based on a different variance estimator,
was found to be no more than 0.63 in that paper).
For design 2 and 3, respectively,
we consider the exchangeable design and the equicorrelated design
studied in Sections~6.1 and 6.2 of \citet{Ber13a}.
The exchangeable design is such that the corresponding PoSI-constant
is small asymptotically,
and the equicorrelated design corresponds to a large PoSI-constant
asymptotically;
cf. Theorem~6.1 and Theorem~6.2 in \citet{Ber13a}.
For the equicorrelated design (design 3), the difference between
the PoSI-interval and the ``naive'' interval is thus expected to be
most pronounced.

More precisely, for the first design, we take the regressor
matrix from a data set from \citet{Raw98a} (page 179), where
the response is peak flow rate
from watersheds, and where the explanatory variables are
rainfall (inches), which is the regressor of interest here,
that is, the first column of $X$, as well as
area of watershed (square miles), area impervious to water (square miles),
average slope of watershed (percent),
longest stream flow in watershed (thousands of feet),
surface absorbency index (0 $=$ complete absorbency; 100 $=$ no absorbency),
estimated soil storage capacity (inches of water),
infiltration rate of water into soil (inches/hour),
time period during which rainfall exceeded $1/4$ inch/hour,
and a constant term to include an intercept in the model.
Logarithms are taken of the response and of all explanatory variables
except for the intercept.
For the second design, we define $\mathbf{X}^{(p)}(a)$ as in Section~6.1
in \citet{Ber13a} with $p=10$ and we choose $a=10$ here, and we
set $X = U \mathbf{X}^{(p)}(a)$, where $U$ is a collection of $p$ orthonormal
$n$-vectors obtained by first drawing a set of $p$ i.i.d. standard
Gaussian $n$-vectors and then applying the Gram--Schmidt procedure.
And for the third design, we define $\mathbf{X}^{(p)}(c)$ as
in Section~6.2 in \citet{Ber13a}, but such that the regressor of
interest is the first one,
where we choose $c = \sqrt{0.8/(p-1)}$, and we set $X = V {\mathbf
X}^{(p)}(c)$,
where $V$ is obtained by drawing an independent observation
from the same distribution as $U$ before.
(Because we consider only orthogonally invariant methods here,
the coverage probabilities under study are invariant under
orthogonal transformations of the columns of the design matrix.
In particular, the coverage probabilities for the second and for
the third design actually do not depend on the matrices $U$ and $V$.)

For each of the three design matrices,
we simulate coverage probabilities under the model $Y=X\beta+ u$ for randomly
selected values of the parameter $\beta$,
we identify those $\beta$'s for which the simulated coverage probability
gets small, and we correct for bias as explained in detail shortly.
For example, consider the case where the coverage target is
$\beta_1$ and where the ``naive'' confidence interval is used with
AIC as the model selector.
We first select 10{,}000 parameters $\beta$
by drawing i.i.d. samples from a random $p$-vector $b$ such
that $X b$ follows a standard Gaussian distribution
within the column-space of $X$.
For each of these $\beta$'s, we approximate the
corresponding coverage probability by the coverage rate obtained
from 100 Monte Carlo samples.
In particular, we draw 100 Monte Carlo samples from the
overall model using $\beta$ as
the true parameter.
For each Monte Carlo sample, we compute the
model selector $\hat{M}$ and the resulting ``naive'' confidence interval,
and we record whether $\beta_1$ is covered or not.
The 100 recorded results are then averaged, resulting in a
coverage rate that provides an
estimator for the coverage probability of the interval if the
true parameter is $\beta$.
After repeating this for each of the 10{,}000 $\beta$'s, we compute
the resulting smallest coverage rate as an estimator for the
minimal coverage probability of the confidence interval.
The smallest coverage rate, as an estimator for the
smallest coverage probability (over the 10,000 selected $\beta$'s),
is
clearly biased downward. To correct for that, we then take those 1000
parameters $\beta$ that gave the smallest coverage rates and re-estimate
the corresponding coverage probabilities as explained earlier, but now
using 1000 Monte Carlo samples.
For that parameter $\beta$ that gives the smallest coverage rate in this
second run, we run the simulation again but now with
500{,}000 Monte Carlo samples, to get a reliable estimate of
the corresponding coverage probability.
This procedure is also used, mutatis mutandis, to evaluate the
performance of
the PoSI1-interval and of the PoSI-interval (with constant $K_{P'}$),
with AIC, BIC, and the LASSO as model selectors,
and also in the case where the coverage target is $\beta_{1\cdot\hat{M}}$.
We stress here that the smallest coverage rates found by this procedure
are simulation-based results obtained by a stochastic search over a
10-dimensional parameter space,
and thus only provide approximate upper bounds for
the true minimal coverage probabilities (cf., e.g., the results
for the PoSI-interval and the PoSI1-interval,
when the coverage target is $\beta_1$,
when BIC is used for model selection, and
when the second design matrix is used for $X$).
Table~\ref{tab1} summarizes the results.
\begin{table*}[t]
\tablewidth=330pt
\caption{Smallest coverage probabilities (rounded to two digits of accuracy
after the comma) found in MC study
for the coverage targets\vspace*{1pt} $\beta_{1\cdot\hat{M}}$,
and $\beta_{1}$,
using AIC, BIC, and the LASSO for model selection,
for the PoSI-interval, the PoSI1-interval, and
the ``naive'' interval, each~with nominal coverage probability 0.95}
\label{tab1}
\begin{tabular*}{330pt}{@{\extracolsep{\fill}}lccccc@{}}
\hline
\textbf{Coverage} & \textbf{Model}
& \textbf{Confidence}
& \textbf{Design 1} & \textbf{Design 2} & \textbf{Design 3}
\\
\textbf{target} & \textbf{selector} & \textbf{interval} & \textbf{(watershed)} & \textbf{(exchangeable)} &
\textbf{(equicorr.)}
\\
\hline
$\beta_{1\cdot\hat{M}}$
&AIC
&PoSI & 1.00 & 1.00 & 0.99
\\
&&PoSI1 & 0.99 & 0.99 & 0.98
\\
&&Naive & 0.89 & 0.92 & 0.81
\\[3pt]
&BIC
&PoSI & 1.00 & 1.00 & 0.99
\\
&&PoSI1 & 0.98 & 0.99 & 0.98
\\
&&Naive & 0.89 & 0.86 & 0.84
\\[3pt]
&LASSO
&PoSI & 1.00 & 1.00 & 1.00
\\
&&PoSI1 & 1.00 & 1.00 & 1.00
\\
&&Naive & 0.95 & 0.95 & 0.93
\\[6pt]
$\beta_{1}$  &AIC
& PoSI & 0.85 & 0.91 & 0.83
\\
&&PoSI1 & 0.76 & 0.91 & 0.77
\\
&& Naive & 0.62 & 0.82 & 0.54
\\[3pt]
&BIC
&PoSI & 0.62 & 0.65 & 0.48
\\
&&PoSI1 & 0.51 & 0.66 & 0.43
\\
&&Naive & 0.43 & 0.51 & 0.26
\\[3pt]
&LASSO
&PoSI & 0.09 & 0.12 & 0.05
\\
&&PoSI1 & 0.08 & 0.12 & 0.03
\\
&&Naive & 0.07 & 0.10 & 0.01
\\
\hline
\end{tabular*}
\end{table*}

For AIC and BIC, the results of the simulation study reinforce the impression
already gained in the theoretical analysis in Section~\ref{theory}:
When the coverage target is $\beta_{1\cdot\hat{M}}$,
the PoSI1-interval as well as the PoSI-interval are somewhat too long and
the ``naive'' interval is somewhat too short, resulting in moderate
over- and under-coverage, respectively.
Both over- and under-coverage are more pronounced than
in the simple model studied in Section~\ref{theory}.
In contrast, when the coverage target is $\beta_1$, then the actual coverage
probability of all intervals can again be far below the nominal level.
As expected, the difference between the ``naive'' interval and
the PoSI1-interval (resp., PoSI-interval) is most pronounced for design 3.
The results for BIC are quite similar to those for AIC, when the
coverage target is $\beta_{1\cdot\hat{M}}$; but when the target
is $\beta_1$, all intervals based on BIC have poorer coverage properties
compared to the intervals based on AIC,
with minima close to, or below, 0.5 in some cases.
This is because BIC selects smaller models than AIC, typically causing more
bias in the resulting post-model-selection estimator
[that phenomenon is analyzed
in greater detail in \citet{Lee03a} and \citet{Poe09a}].
The results for the LASSO stand out: When the coverage target is
$\beta_{1\cdot\hat{M}}$, the PoSI1-interval (resp., PoSI-interval)
gives smallest probabilities very close to one, while the smallest coverage
probability of the naive interval is very close to the nominal level (0.95).
But when the coverage target is $\beta_{1}$, all intervals have
smallest coverage probabilities of around $0.1$ and below.
The reason for this is that the LASSO model selector, as implemented here
and for the parameters used in the stochastic search for the smallest
coverage probability, selects the smallest possible model in most cases,
that is, the model containing only the first regressor.
In other words, the model selected by the LASSO is ``nearly nonrandom.''
When the target is $\beta_{1\cdot\hat{M}}$, this entails that the
naive interval is approximately valid and that both PoSI intervals are
too large. [Indeed, the naive interval is valid if the underlying model
selector always chooses a fixed (nonrandom) model;
cf. the discussion following \eqref{CI}.]
But when the target is $\beta_1$, the model selected by the LASSO typically
suffers from severe bias, resulting in very small coverage probabilities
for all intervals.

Other model selectors can, of course, give results different from
those in Table~\ref{tab1}. The model selectors chosen here represent
a selection of popular methods from the contemporary literature that
exhibit an interesting range of possible scenarios for the
minimal coverage probabilities of confidence intervals post-model-selection.

\begin{appendix}

\section*{Appendix: Confidence Sets Under Zero-Restrictions
Post-Model-Selection}\label{app}

Let $y$ and $\hat{\sigma}^2$ be as in Section~\ref{setting},
and consider
$\mathcal{M} = \{M_0, M_1\}$, where each of the two candidate
models $M_i$ is full rank.
Suppose we are interested in the coefficient of the first regressor
$X_1$, that is assumed present in $M_1$ but absent in $M_0$.
In the notation introduced in Section~\ref{setting}, we thus have
$1 \in M_1$ and $1 \notin M_0$.
Let $\hat{M}$ be \textit{any} model selection procedure
that chooses only between $M_0$ and $M_1$.
As the model-dependent coverage target, we consider the
coefficient of $X_1$, which is
not restricted under $M_1$, and which is restricted to zero under $M_0$.
More precisely, set $b_{M_1} = \beta_{1\cdot M_1}$,
set $b_{M_0} = 0$, and let the target be $b_{\hat{M}}$.
We consider a ``naive'' confidence interval for $b_{\hat{M}}$ that is
defined as
\[
I_{\hat{M}} = \cases{\ds
\hat{\beta}_{1\cdot M_1}
\pm k_N \hat{\sigma}_{1\cdot M_1},  & \mbox{if $\hat{M} =
M_1$},
\vspace*{2pt}\cr
\ds \{0\},  & \mbox{if $\hat{M}=M_0$,}}
\]
where $k_N$ is chosen so that
$\mathbb{P}_{\mu,\sigma}( \beta_{1\cdot M_1} \in \hat{\beta
}_{1\cdot M_1}
\pm k_N \hat{\sigma}_{1\cdot M_1}) = 1-\alpha$.
[The constant $k_N$ is
the $(1-\alpha/2)$-quantile of
a standard normal distribution in the known-variance
case and
the $(1-\alpha/2)$-quantile of
a $t$-distribution with $r$ degrees of freedom in the unknown-variance case.]
The actual coverage probability of $I_{\hat{M}}$, as a confidence interval
for $b_{\hat{M}}$, is at least equal to the nominal coverage
probability $1-\alpha$, because
\begin{eqnarray*}
&& \mathbb{P}_{\mu,\sigma}( b_{\hat{M}} \in I_{\hat{M}})
\\
&&\quad = \mathbb{P}_{\mu,\sigma} (\beta_{1\cdot M_1} \in I_{M_1}\mbox{
and } \hat{M}=M_1) \\
&&\qquad{}+ \mathbb{P}_{\mu,\sigma} \bigl( 0 \in\{0\},
\hat{M}=M_0 \bigr)
\\
&&\quad = \mathbb{P}_{\mu,\sigma} (\beta_{1\cdot M_1} \in I_{M_1}\mbox{
and } \hat{M}= M_1)\\
&& \qquad{}+ \mathbb{P}_{\mu,\sigma}( \hat{M}\neq
M_1)
\\
&&\quad = \mathbb{P}_{\mu,\sigma}( \beta_{1\cdot M_1} \in I_{M_1} \mbox{
or } \hat{M}\neq M_1 ) \geq 1-\alpha,
\end{eqnarray*}
where the inequality in the last step
holds in view of the choice of $k_N$.
\end{appendix}

\section*{Acknowledgments}

We thank the anonymous referee and the Editor for helpful comments
and feedback. We also thank
Richard Berk, Lawrence Brown, Andreas Buja, Kai Zhang, and Linda Zhao
for providing us with the code to compute the
PoSI-constant $K_{P'}$ used in Section~\ref{simulation};
the entire ``PoSI-group'' at the University of Pennsylvania for
inspiring discussions during Hannes Leeb's visit;
and Francois Bachoc for constructive feedback.

Karl Ewald supported in part by Deutsche\break Forschungsgemeinschaft (DFG)
Grant FOR916, and
Hannes Leeb supported in part by FWF Grant P26354.



%



\begin{thebibliography}{32}

\bibitem[\protect\citeauthoryear{Andrews and Guggenberger}{2009}]{And09a}
\begin{barticle}[mr]
\bauthor{\bsnm{Andrews},~\bfnm{Donald~W.~K.}\binits{D.~W.~K.}} \AND
\bauthor{\bsnm{Guggenberger},~\bfnm{Patrik}\binits{P.}}
(\byear{2009}).
\btitle{Hybrid and size-corrected subsampling methods}.
\bjournal{Econometrica}
\bvolume{77}
\bpages{721--762}.
\bid{doi={10.3982/ECTA7015}, issn={0012-9682}, mr={2531360}}
\end{barticle}
%

\bptok{imsref}%
\endbibitem

\bibitem[\protect\citeauthoryear{Berk et~al.}{2013}]{Ber13a}
\begin{barticle}[mr]
\bauthor{\bsnm{Berk},~\bfnm{Richard}\binits{R.}},
\bauthor{\bsnm{Brown},~\bfnm{Lawrence}\binits{L.}},
\bauthor{\bsnm{Buja},~\bfnm{Andreas}\binits{A.}},
\bauthor{\bsnm{Zhang},~\bfnm{Kai}\binits{K.}} \AND
\bauthor{\bsnm{Zhao},~\bfnm{Linda}\binits{L.}}
(\byear{2013}).
\btitle{Valid post-selection inference}.
\bjournal{Ann. Statist.}
\bvolume{41}
\bpages{802--837}.
\bid{doi={10.1214/12-AOS1077}, issn={0090-5364}, mr={3099122}}
\end{barticle}
%

\bptok{imsref}%
\endbibitem

\bibitem[\protect\citeauthoryear{Bickel and Doksum}{1977}]{Bic77a}
\begin{bbook}[mr]
\bauthor{\bsnm{Bickel},~\bfnm{Peter~J.}\binits{P.~J.}} \AND
\bauthor{\bsnm{Doksum},~\bfnm{Kjell~A.}\binits{K.~A.}}
(\byear{1977}).
\btitle{Mathematical Statistics: Basic Ideas and Selected Topics}.
\bpublisher{Holden-Day},
\blocation{Oakland, CA}.
\bid{mr={0443141}}
\bptnote{check year}%
\end{bbook}
%

\bptok{imsref}%
\endbibitem

\bibitem[\protect\citeauthoryear{Brown}{1967}]{Bro67a}
\begin{barticle}[mr]
\bauthor{\bsnm{Brown},~\bfnm{L.}\binits{L.}}
(\byear{1967}).
\btitle{The conditional level of {S}tudent's {$t$} test}.
\bjournal{Ann. Math. Stat.}
\bvolume{38}
\bpages{1068--1071}.
\bid{issn={0003-4851}, mr={0214210}}
\end{barticle}
%

\bptok{imsref}%
\endbibitem

\bibitem[\protect\citeauthoryear{Buehler and Feddersen}{1963}]{Bue63a}
\begin{barticle}[mr]
\bauthor{\bsnm{Buehler},~\bfnm{R.~J.}\binits{R.~J.}} \AND
\bauthor{\bsnm{Feddersen},~\bfnm{A.~P.}\binits{A.~P.}}
(\byear{1963}).
\btitle{Note on a conditional property of {S}tudent's {$t$}}.
\bjournal{Ann. Math. Stat.}
\bvolume{34}
\bpages{1098--1100}.
\bid{issn={0003-4851}, mr={0150864}}
\end{barticle}
%

\bptok{imsref}%
\endbibitem

\bibitem[\protect\citeauthoryear{Craven and Wahba}{1978/79}]{Cra78a}
\begin{barticle}[mr]
\bauthor{\bsnm{Craven},~\bfnm{Peter}\binits{P.}} \AND
\bauthor{\bsnm{Wahba},~\bfnm{Grace}\binits{G.}}
(\byear{1978/79}).
\btitle{Smoothing noisy data with spline functions. {E}stimating the correct degree of smoothing by the method of generalized cross-validation}.
\bjournal{Numer. Math.}
\bvolume{31}
\bpages{377--403}.
\bid{doi={10.1007/BF01404567}, issn={0029-599X}, mr={0516581}}
\bptnote{check year}%
\end{barticle}
%

\bptok{imsref}%
\endbibitem

\bibitem[\protect\citeauthoryear{Dijkstra and Veldkamp}{1988}]{Dij88a}
\begin{bincollection}[auto:parserefs-M02]
\bauthor{\bsnm{Dijkstra},~\bfnm{T.~K.}\binits{T.~K.}} \AND
\bauthor{\bsnm{Veldkamp},~\bfnm{J.~H.}\binits{J.~H.}}
(\byear{1988}).
\btitle{Data-driven selection of regressors and the bootstrap}.
In \bbooktitle{Lecture Notes in Econom. and Math. Systems}
\bvolume{307}
\bpages{17--38}.
\bpublisher{Springer},
\blocation{New York}.
\end{bincollection}
%

\bptok{imsref}%
\endbibitem

\bibitem[\protect\citeauthoryear{Efron et~al.}{2004}]{Efr04a}
\begin{barticle}[mr]
\bauthor{\bsnm{Efron},~\bfnm{Bradley}\binits{B.}},
\bauthor{\bsnm{Hastie},~\bfnm{Trevor}\binits{T.}},
\bauthor{\bsnm{Johnstone},~\bfnm{Iain}\binits{I.}} \AND
\bauthor{\bsnm{Tibshirani},~\bfnm{Robert}\binits{R.}}
(\byear{2004}).
\btitle{Least angle regression}.
\bjournal{Ann. Statist.}
\bvolume{32}
\bpages{407--499}.
\bid{doi={10.1214/009053604000000067}, issn={0090-5364}, mr={2060166}}
\bptnote{check related}%
\end{barticle}
%

\bptok{imsref}%
\endbibitem

\bibitem[\protect\citeauthoryear{Ewald}{2012}]{Ewa12a}
\begin{bmisc}[auto]
\bauthor{\bsnm{Ewald},~\bfnm{K.}\binits{K.}}
(\byear{2012}).
\bhowpublished{On the influence of model selection on confidence regions for marginal associations in the linear
model.
Master's thesis,
Univ. Vienna.}
\end{bmisc}
%

\bptok{imsref}%
\endbibitem

\bibitem[\protect\citeauthoryear{Kabaila}{1998}]{Kab98a}
\begin{barticle}[mr]
\bauthor{\bsnm{Kabaila},~\bfnm{Paul}\binits{P.}}
(\byear{1998}).
\btitle{Valid confidence intervals in regression after variable selection}.
\bjournal{Econometric Theory}
\bvolume{14}
\bpages{463--482}.
\bid{doi={10.1017/S0266466698144031}, issn={0266-4666}, mr={1650037}}
\end{barticle}
%

\bptok{imsref}%
\endbibitem

\bibitem[\protect\citeauthoryear{Kabaila}{2009}]{Kab09a}
\begin{barticle}[auto:parserefs-M02]
\bauthor{\bsnm{Kabaila},~\bfnm{P.}\binits{P.}}
(\byear{2009}).
\btitle{The coverage properties of confidence regions after model selection}.
\bjournal{Int. Stat. Rev.}
\bvolume{77}
\bpages{405--414}.
\end{barticle}
%

\bptok{imsref}%
\endbibitem

\bibitem[\protect\citeauthoryear{Kabaila and Leeb}{2006}]{Kab06a}
\begin{barticle}[mr]
\bauthor{\bsnm{Kabaila},~\bfnm{Paul}\binits{P.}} \AND
\bauthor{\bsnm{Leeb},~\bfnm{Hannes}\binits{H.}}
(\byear{2006}).
\btitle{On the large-sample minimal coverage probability of confidence intervals after model selection}.
\bjournal{J. Amer. Statist. Assoc.}
\bvolume{101}
\bpages{619--629}.
\bid{doi={10.1198/016214505000001140}, issn={0162-1459}, mr={2256178}}
\end{barticle}
%

\bptok{imsref}%
\endbibitem

\bibitem[\protect\citeauthoryear{Leeb}{2006}]{Lee01d}
\begin{bincollection}[mr]
\bauthor{\bsnm{Leeb},~\bfnm{Hannes}\binits{H.}}
(\byear{2006}).
\btitle{The distribution of a linear predictor after model selection: Unconditional finite-sample distributions and asymptotic approximations}.
In \bbooktitle{Optimality}.
\bseries{Institute of Mathematical Statistics Lecture Notes---Monograph Series}
\bvolume{49}
\bpages{291--311}.
\bpublisher{IMS},
\blocation{Beachwood, OH}.
\bid{doi={10.1214/074921706000000518}, mr={2338549}}
\end{bincollection}
%

\bptok{imsref}%
\endbibitem

\bibitem[\protect\citeauthoryear{Leeb}{2008}]{Lee05a}
\begin{barticle}[mr]
\bauthor{\bsnm{Leeb},~\bfnm{Hannes}\binits{H.}}
(\byear{2008}).
\btitle{Evaluation and selection of models for out-of-sample prediction when the sample size is small relative to the complexity of the data-generating process}.
\bjournal{Bernoulli}
\bvolume{14}
\bpages{661--690}.
\bid{doi={10.3150/08-BEJ127}, issn={1350-7265}, mr={2537807}}
\end{barticle}
%

\bptok{imsref}%
\endbibitem

\bibitem[\protect\citeauthoryear{Leeb and P{\"o}tscher}{2003}]{Lee02d}
\begin{barticle}[mr]
\bauthor{\bsnm{Leeb},~\bfnm{Hannes}\binits{H.}} \AND
\bauthor{\bsnm{P{\"o}tscher},~\bfnm{Benedikt~M.}\binits{B.~M.}}
(\byear{2003}).
\btitle{The finite-sample distribution of post-model-selection estimators and uniform versus nonuniform approximations}.
\bjournal{Econometric Theory}
\bvolume{19}
\bpages{100--142}.
\bid{doi={10.1017/S0266466603191050}, issn={0266-4666}, mr={1965844}}
\end{barticle}
%

\bptok{imsref}%
\endbibitem

\bibitem[\protect\citeauthoryear{Leeb and P{\"o}tscher}{2005}]{Lee03a}
\begin{barticle}[mr]
\bauthor{\bsnm{Leeb},~\bfnm{Hannes}\binits{H.}} \AND
\bauthor{\bsnm{P{\"o}tscher},~\bfnm{Benedikt~M.}\binits{B.~M.}}
(\byear{2005}).
\btitle{Model selection and inference: Facts and fiction}.
\bjournal{Econometric Theory}
\bvolume{21}
\bpages{21--59}.
\bid{doi={10.1017/S0266466605050036}, issn={0266-4666}, mr={2153856}}
\end{barticle}
%

\bptok{imsref}%
\endbibitem

\bibitem[\protect\citeauthoryear{Leeb and P{\"o}tscher}{2006a}]{Lee02a}
\begin{barticle}[mr]
\bauthor{\bsnm{Leeb},~\bfnm{Hannes}\binits{H.}} \AND
\bauthor{\bsnm{P{\"o}tscher},~\bfnm{Benedikt~M.}\binits{B.~M.}}
(\byear{2006}a).
\btitle{Can one estimate the conditional distribution of post-model-selection estimators?}
\bjournal{Ann. Statist.}
\bvolume{34}
\bpages{2554--2591}.
\bid{doi={10.1214/009053606000000821}, issn={0090-5364}, mr={2291510}}
\end{barticle}
%

\bptok{imsref}%
\endbibitem

\bibitem[\protect\citeauthoryear{Leeb and P{\"o}tscher}{2006b}]{Lee02c}
\begin{barticle}[mr]
\bauthor{\bsnm{Leeb},~\bfnm{Hannes}\binits{H.}} \AND
\bauthor{\bsnm{P{\"o}tscher},~\bfnm{Benedikt~M.}\binits{B.~M.}}
(\byear{2006}b).
\btitle{Performance limits for estimators of the risk or distribution of shrinkage-type estimators, and some general lower risk-bound results}.
\bjournal{Econometric Theory}
\bvolume{22}
\bpages{69--97}.
\bid{doi={10.1017/S0266466606060038}, issn={0266-4666}, mr={2212693}}
\end{barticle}
%

\bptok{imsref}%
\endbibitem

\bibitem[\protect\citeauthoryear{Leeb and P{\"o}tscher}{2008a}]{Lee02b}
\begin{barticle}[mr]
\bauthor{\bsnm{Leeb},~\bfnm{Hannes}\binits{H.}} \AND
\bauthor{\bsnm{P{\"o}tscher},~\bfnm{Benedikt~M.}\binits{B.~M.}}
(\byear{2008}a).
\btitle{Can one estimate the unconditional distribution of post-model-selection estimators?}
\bjournal{Econometric Theory}
\bvolume{24}
\bpages{338--376}.
\bid{doi={10.1017/S0266466608080158}, issn={0266-4666}, mr={2422862}}
\end{barticle}
%

\bptok{imsref}%
\endbibitem

\bibitem[\protect\citeauthoryear{Leeb and P{\"{o}}tscher}{2008b}]{Lee06a}
\begin{bincollection}[auto:parserefs-M02]
\bauthor{\bsnm{Leeb},~\bfnm{H.}\binits{H.}} \AND
\bauthor{\bsnm{P{\"{o}}tscher},~\bfnm{B.~M.}\binits{B.~M.}}
(\byear{2008}b).
\btitle{Model selection}.
In \bbooktitle{Handbook of Financial Time Series}
(\beditor{\bfnm{T.~G.}\binits{T.~G.}~\bsnm{Andersen}},
\beditor{\bfnm{R.~A.}\binits{R.~A.}~\bsnm{Davis}},
\beditor{\bfnm{J.-P.}\binits{J.-P.}~\bsnm{Krei{\ss}}} \AND
\beditor{\bfnm{Th.}\binits{Th.}~\bsnm{Mikosch}}, eds.)
\bpages{785--821}.
\bpublisher{Springer},
\blocation{New York}.
\end{bincollection}
%

\bptok{imsref}%
\endbibitem

\bibitem[\protect\citeauthoryear{Olshen}{1973}]{Ols73a}
\begin{barticle}[mr]
\bauthor{\bsnm{Olshen},~\bfnm{Richard~A.}\binits{R.~A.}}
(\byear{1973}).
\btitle{The conditional level of the {$F$}-test}.
\bjournal{J.~Amer. Statist. Assoc.}
\bvolume{68}
\bpages{692--698}.
\bid{issn={0162-1459}, mr={0359198}}
\end{barticle}
%

\bptok{imsref}%
\endbibitem

\bibitem[\protect\citeauthoryear{P{\"o}tscher}{1991}]{Poe91a}
\begin{barticle}[mr]
\bauthor{\bsnm{P{\"o}tscher},~\bfnm{B.~M.}\binits{B.~M.}}
(\byear{1991}).
\btitle{Effects of model selection on inference}.
\bjournal{Econometric Theory}
\bvolume{7}
\bpages{163--185}.
\bid{doi={10.1017/S0266466600004382}, issn={0266-4666}, mr={1128410}}
\end{barticle}
%

\bptok{imsref}%
\endbibitem

\bibitem[\protect\citeauthoryear{P{\"o}tscher}{2006}]{Poe06a}
\begin{bincollection}[mr]
\bauthor{\bsnm{P{\"o}tscher},~\bfnm{Benedikt~M.}\binits{B.~M.}}
(\byear{2006}).
\btitle{The distribution of model averaging estimators and an impossibility result regarding its estimation}.
In \bbooktitle{Time Series and Related Topics}.
\bseries{Institute of Mathematical Statistics Lecture Notes---Monograph Series}
\bvolume{52}
\bpages{113--129}.
\bpublisher{IMS},
\blocation{Beachwood, OH}.
\bid{doi={10.1214/074921706000000987}, mr={2427842}}
\end{bincollection}
%

\bptok{imsref}%
\endbibitem

\bibitem[\protect\citeauthoryear{P{\"o}tscher}{2009}]{Poe09a}
\begin{barticle}[mr]
\bauthor{\bsnm{P{\"o}tscher},~\bfnm{Benedikt~M.}\binits{B.~M.}}
(\byear{2009}).
\btitle{Confidence sets based on sparse estimators are necessarily large}.
\bjournal{Sankhy\={a}}
\bvolume{71}
\bpages{1--18}.
\bid{issn={0972-7671}, mr={2579644}}
\end{barticle}
%

\bptok{imsref}%
\endbibitem

\bibitem[\protect\citeauthoryear{P{\"o}tscher and Leeb}{2009}]{Poe09b}
\begin{barticle}[mr]
\bauthor{\bsnm{P{\"o}tscher},~\bfnm{Benedikt~M.}\binits{B.~M.}} \AND
\bauthor{\bsnm{Leeb},~\bfnm{Hannes}\binits{H.}}
(\byear{2009}).
\btitle{On the distribution of penalized maximum likelihood estimators: The LASSO, SCAD, and thresholding}.
\bjournal{J. Multivariate Anal.}
\bvolume{100}
\bpages{2065--2082}.
\bid{doi={10.1016/j.jmva.2009.06.010}, issn={0047-259X}, mr={2543087}}
\end{barticle}
%

\bptok{imsref}%
\endbibitem

\bibitem[\protect\citeauthoryear{P{\"o}tscher and Schneider}{2009}]{Poe09c}
\begin{barticle}[mr]
\bauthor{\bsnm{P{\"o}tscher},~\bfnm{Benedikt~M.}\binits{B.~M.}} \AND
\bauthor{\bsnm{Schneider},~\bfnm{Ulrike}\binits{U.}}
(\byear{2009}).
\btitle{On the distribution of the adaptive LASSO estimator}.
\bjournal{J. Statist. Plann. Inference}
\bvolume{139}
\bpages{2775--2790}.
\bid{doi={10.1016/j.jspi.2009.01.003}, issn={0378-3758}, mr={2523666}}
\end{barticle}
%

\bptok{imsref}%
\endbibitem

\bibitem[\protect\citeauthoryear{P{\"o}tscher and Schneider}{2010}]{Poe10a}
\begin{barticle}[mr]
\bauthor{\bsnm{P{\"o}tscher},~\bfnm{Benedikt~M.}\binits{B.~M.}} \AND
\bauthor{\bsnm{Schneider},~\bfnm{Ulrike}\binits{U.}}
(\byear{2010}).
\btitle{Confidence sets based on penalized maximum likelihood estimators in {G}aussian regression}.
\bjournal{Electron. J. Stat.}
\bvolume{4}
\bpages{334--360}.
\bid{doi={10.1214/09-EJS523}, issn={1935-7524}, mr={2645488}}
\end{barticle}
%

\bptok{imsref}%
\endbibitem

\bibitem[\protect\citeauthoryear{P{\"o}tscher and Schneider}{2011}]{Poe11a}
\begin{barticle}[mr]
\bauthor{\bsnm{P{\"o}tscher},~\bfnm{Benedikt~M.}\binits{B.~M.}} \AND
\bauthor{\bsnm{Schneider},~\bfnm{Ulrike}\binits{U.}}
(\byear{2011}).
\btitle{Distributional results for thresholding estimators in high-dimensional {G}aussian regression models}.
\bjournal{Electron. J. Stat.}
\bvolume{5}
\bpages{1876--1934}.
\bid{doi={10.1214/11-EJS659}, issn={1935-7524}, mr={2970179}}
\end{barticle}
%

\bptok{imsref}%
\endbibitem

\bibitem[\protect\citeauthoryear{Rawlings, Pantula and Dickey}{1998}]{Raw98a}
\begin{bbook}[mr]
\bauthor{\bsnm{Rawlings},~\bfnm{John~O.}\binits{J.~O.}},
\bauthor{\bsnm{Pantula},~\bfnm{Sastry~G.}\binits{S.~G.}} \AND
\bauthor{\bsnm{Dickey},~\bfnm{David~A.}\binits{D.~A.}}
(\byear{1998}).
\btitle{Applied Regression Analysis: A Research Tool},
\bedition{2nd} ed.
\bpublisher{Springer},
\blocation{New York}.
\bid{doi={10.1007/b98890}, mr={1631919}}
\end{bbook}
%

\bptok{imsref}%
\endbibitem

\bibitem[\protect\citeauthoryear{Sen}{1979}]{Sen79a}
\begin{barticle}[mr]
\bauthor{\bsnm{Sen},~\bfnm{Pranab~Kumar}\binits{P.~K.}}
(\byear{1979}).
\btitle{Asymptotic properties of maximum likelihood estimators based on conditional specification}.
\bjournal{Ann. Statist.}
\bvolume{7}
\bpages{1019--1033}.
\bid{issn={0090-5364}, mr={0536504}}
\end{barticle}
%

\bptok{imsref}%
\endbibitem

\bibitem[\protect\citeauthoryear{Sen and Saleh}{1987}]{Sen87a}
\begin{barticle}[mr]
\bauthor{\bsnm{Sen},~\bfnm{Pranab~Kumar}\binits{P.~K.}} \AND
\bauthor{\bsnm{Saleh},~\bfnm{A.~K.~M.~Ehsanes}\binits{A.~K.~M.~E.}}
(\byear{1987}).
\btitle{On preliminary test and shrinkage {$M$}-estimation in linear models}.
\bjournal{Ann. Statist.}
\bvolume{15}
\bpages{1580--1592}.
\bid{doi={10.1214/aos/1176350611}, issn={0090-5364}, mr={0913575}}
\end{barticle}
%

\bptok{imsref}%
\endbibitem

\bibitem[\protect\citeauthoryear{Tukey}{1967}]{Tuc67a}
\begin{barticle}[auto]
\bauthor{\bsnm{Tukey},~\bfnm{J.~W.}\binits{J.~W.}}
(\byear{1967}).
\btitle{Discussion of ``Topics in the investigation of linear relations
fitted by the method of least squares'' by F.~J.~Anscombe}.
\bjournal{J. Roy. Statist. Soc. Ser. B}
\bvolume{29}
\bpages{47--48}.
\end{barticle}
%
\bptok{imsref}%
\endbibitem

\end{thebibliography}
\end{document}